\documentclass[12pt,reqno]{amsart}
\usepackage{epsf, amsmath, amsfonts, amsthm, amssymb}
\theoremstyle{plain}
\newtheorem{theorem}{Theorem}[section]
\newtheorem{lemma}[theorem]{Lemma}
\newtheorem{proposition}[theorem]{Proposition}

\newtheorem{definition}[theorem]{Definition}
                                                                                
\theoremstyle{definition}
\newtheorem{observation}[theorem]{Observation}
\newtheorem{remark}[theorem]{Remark}

\begin{document}

\newcommand{\gruppo}{\mathbb{Z}_3}
\newcommand{\luna}{\mathbb{Z}_9}
\newcommand{\freefactor}{\mathfrak{L} \left (\mathbf{F} _\frac{11}{3}\right)}
\newcommand{\inff}{\mathfrak{L} \left (\mathbf{F} _{t}\right)}
\newcommand{\pippo}{\inff\otimes R}
\newcommand{\tenpro}{\freefactor\otimes R_{0}}
\newcommand{\pten}{\left (\tenpro\right )}
\newcommand{\crossprod}{\left ( \tenpro\right )\rtimes_{\gamma}\gruppo}
\newcommand{\obstruction}{e^{\frac{2\pi i}{3}}} 
\newcommand{\obstruconj}{e^{\frac{-2\pi i}{3}}} 
\newcommand{\jonesinv}{e^{\frac{2\pi i}{9}}}
\newcommand{\jonesconj}{e^{\frac{-2\pi i}{9}}}
\newcommand{\basic}{R_{-1}\rtimes _{\theta}\luna}
\newcommand{\matrices}{M_{9}(\mathbb{C})}
\newcommand{\jconj}{\bar{\delta}}
\newcommand{\M}{\mathcal {M}}
\newcommand{\Ad}{\operatorname{Ad\,}}
\newcommand{\fint}{\{\ad _{N}u\, |\, u\in N\text{ is fixed by }G\}}
\newcommand{\closure}{\overline{\mbox{Fint}}}
\newcommand{\clos}{clos\{\Ad u\, |\, u \text{ is fixed by }\gruppo \}}
\newcommand{\gener}{Id\otimes (\Ad U_0 ^{*}\,\beta)}
\newcommand{\dualgr}{\widehat{\mathbb{Z}}_{3}}
\newcommand{\dualcr}{\M\rtimes_{\widehat{\gamma}}\dualgr} 
\newcommand{\qvd}{$\hfill\blacksquare$}
\newcommand{\ad}{\operatorname{Ad}}
\newcommand{\Rea}{\mathrm{Re}}
%
%
\newcommand{\innerM}{\operatorname{Int}(M)}
\newcommand{\innerN}{\operatorname{Int}(N)}
\newcommand{\innerR}{\operatorname{Int}(R)}
\newcommand{\innerfac}{\operatorname{Int}\left (\freefactor\right )}
\newcommand{\innerpten}{\operatorname{Int}\pten}
\newcommand{\CtM}{\operatorname{Ct}(M)}
\newcommand{\CtN}{\operatorname{Ct}(N)}
\newcommand{\CtR}{\operatorname{Ct}(R_{0})}
\newcommand{\Ctpten}{\operatorname{Ct}\pten}
\newcommand{\Ctpi}{\operatorname{Ct}\left (\pippo\right )}
\newcommand{\autM}{\operatorname{Aut}(M)} 
\newcommand{\autR}{\operatorname{Aut}(R)}
\newcommand{\autpten}{\operatorname{Aut}\pten}
\newcommand{\autfac}{\operatorname{Aut}\left (\freefactor\right )}
\newcommand{\autN}{\operatorname{Aut}(N)} 
\newcommand{\outM}{\operatorname{Out}(M)}
\newcommand{\outN}{\operatorname{Out}(N)}
\newcommand{\chiusoM}{\overline{\operatorname {Int}(M)}}
\newcommand{\chiusoN}{\overline{\operatorname{Int}(N)}}
\newcommand{\chiusoR}{\overline{\operatorname{Int}(R)}}
\newcommand{\chiusofac}{\overline{\operatorname{Int}\left (\freefactor\right )}}
\newcommand{\chiusopten}{\overline{\operatorname{Int}\pten}}
%
%
%
\title[On a Subfactor Construction]{On a Subfactor Construction of a Factor Non-Antiisomorphic to Itself}
\author{by Maria Grazia Viola}
\address{\hskip-\parindent Department of Mathematics \\Texas A\&M University \\ College Station TX 77843, USA\\ Fax: (979) 845-3643.} 
\email{viola@math.tamu.edu}

\begin{abstract}
We define a $\gruppo$-kernel  $\alpha$ on $\freefactor$ and a $\gruppo$-kernel $\beta$ on the 
hyperfinite factor $R$, which have conjugate obstruction to lifting. Hence, $\alpha\otimes\beta$ can be perturbed by 
 an inner automorphism to produce an action $\gamma$ on $\tenpro$. The aim of this paper is to show that the factor 
$\M=\crossprod$, which is similar to Connes's example of a $II_{1}$ factor non-antiisomorphic to itself, is the enveloping algebra of an 
inclusion of $II_{1}$ factors $A\subset B$. Here $A$ is a free group 
factor and $B$ is isomorphic to the crossed product $A\rtimes_{\theta}\luna$, 
where $\theta$ is a $\gruppo$-kernel of $A$\ with non-trivial obstruction to 
lifting. By using an argument due to Connes, which involves the invariant $\chi (\M)$, we 
show that $\M$ is not anti-isomorphic to itself. Furthermore, we
prove that for one of the generator of $\chi (\M)$, which we will denote by $\sigma$, 
the Jones invariant $\varkappa (\sigma)$ is equal to $\jonesinv$.
\end{abstract}

\maketitle 

\section{Introduction}

A von Neumann algebra $M$ is anti-isomorphic to itself if there exists a vector space isomorphism $\Phi :M\longrightarrow M$ with the properties $\Phi (x^{*})=\Phi (x)^*$ and $\Phi (xy)=\Phi (y)\Phi (x)$ for every $x,y\in M$. This is equivalent to saying that  $M$ is isomorphic to its conjugate algebra $M^{c}$ (defined in Section 6). 

With his examples of a II$_{1}$ factor non-antiisomorphic to itself (cf. \cite{Connes6}), A. Connes gave in the '70s an answer to a crucial  problem posed by F. Murray and J. von Neumann a few decades earlier, concerning the possibility of realizing every II$_{1}$ factor as the left regular representation of a discrete group. His example was obtained from the II$_{1}$ factor $\mathfrak{L} \left (\mathbf{F} _{4}\right)\otimes R$, where $R$ denotes the hyperfinite II$_{1}$ factor, using a crossed product construction with a $\gruppo$-action. After the innovative work on subfactors done by V. Jones in the '80s (see \cite{Jones3}), it was a natural question to ask whether Connes's factor could be obtained through a subfactor construction of finite index. Although extensive work (\cite{Connes1}, \cite{Connes6}, and \cite{Jones1}) has been done by both Connes and Jones on examples of II$_{1}$ and III$_{\lambda}$ factors non-antiisomorphic to itself, it does not seems that this problem has been addressed before and there is little literature on the subject. 

In this paper we provide a positive answer to this question by giving an explicit subfactor construction for our example of a II$_{1}$ factor non-antiisomorphic to itself. Our model is a variation of Connes's example (\cite{Connes4} and 
\cite{Connes6}), since we utilize in our approach the recently developed theory of interpolated free group factors (\cite{Radulescu1} and \cite{Dykema2}). The II$_{1}$ factor $\M$ we are going to study is constructed from the tensor product of the interpolated free group factor $\freefactor$ and the hyperfinite II$_{1}$ factor $R$. We use two $\gruppo$-kernels, $\alpha\in\autfac$ and $\beta\in\autR$, 
which have conjugate obstructions to lifting, to generate an action of $\gruppo$ on $\freefactor\otimes R$. The action is given, up to 
an inner automorphism, by $\alpha\otimes\beta$, and the factor $\M$ is equal to the crossed product $\left (\freefactor\otimes R\right )\rtimes _{\gamma}\gruppo$ (cf. Section $4$).
     
The main result of this paper is Theorem \ref{main}, where we show that $\M$ is the enveloping algebra of an inclusion  
$A\subset B$ of interpolated free group factors. Here $A$ is isomorphic to $\mathfrak{L}\left (\mathbf{F} _{\frac{35}{27}} \right)$ (Proposition \ref{proposition4.3}), 
and $B$ is equal to the the crossed product $A\rtimes _{\theta}\luna$, for a $\luna$-action $\theta$ of $A$ with outer period 3, and obstruction $\obstruction$ to lifting. The proof is based on Voiculescu's random matrix model for circular and semicircular elements introduced 
in \cite{Vocu}. An analogous argument has been used by F. R\v{a}dulescu in \cite{Radulescu2}, to prove that a variation of the example 
given by Jones of a II$_{1}$ factor with Connes invariant equal to $\mathbb{Z}_{2}\otimes\mathbb{Z} _{2}$, has a subfactor construction. 
We also show in Section 5 that the Connes invariant of our factor $\M$ is equal to $\luna$, a result announced by Connes in \cite{Connes4}. This invariant, which is defined for every factor $M$ with a separable predual, was introduced by Connes in \cite{Connes6}. It consists of an abelian subgroup of the group of outer automorphisms, and it is denoted by $\chi (M)$. It is an important tool for distinguishing factors, since it is an isomorphism invariant of the factor $M$. It is trivial for some of the most common II$_{1}$ factors, like the interpolated free group factors and the hyperfinite II$_{1}$ factor, as well as for any tensor product of these factors. However, a crossed product construction yields in general a non-trivial $\chi (M)$. To 
compute the Connes invariant of the factor $\M=\crossprod$, we use the short exact sequence described by Connes in \cite{Connes6}.  

In addition, we show that if $\sigma$ denotes the generator of $\chi (\M)$ described in Remark \ref{add}, then the invariant $\varkappa (\sigma)$, introduced by Jones in \cite{Jones1} is equal to $\jonesinv$. The Jones invariant is defined for any element $\theta$ of $\chi (M)$, where $M$ is a II$_{1}$ factor without non-trivial hypercentral sequences, and it consists of a complex number of modulus one. It is a finer invariant than $\chi (M)$ and it is constant on the conjugacy class of $\theta$ in the group of outer automorphisms. Moreover, it behaves nicely with respect to antiautomorphisms of $M$, in the sense that conjugation by an anti automorphism changes  $\varkappa (\theta )$ by complex conjugation (see \cite{Jones1} for details).
       
Lastly, in Section 6 we use an argument of Connes \cite{Connes4} to show that $\M=$\linebreak $\crossprod$ is not anti-isomorphic to itself. The two main ingredients of this argument are the uniqueness (up to inner automorphism) of the decomposition of $\gamma$ into the product of an approximately inner automorphism and a centrally trivial automorphism, and the fact that the unique subgroup of order $3$ in $\chi (\M)$ is 
generated by the dual action $\widehat{\gamma}$ on $\M=\crossprod$. Using this decomposition we obtain a canonical way to associate to the II$_{1}$ factor $\M$ a complex number, the obstruction to lifting of the approximately inner automorphism in the decomposition of $\gamma$, which is invariant under isomorphism, and distinguishes $\M$ from its conjugate algebra $\M ^{c}$.

\section{Definitions}

Let $M$ be a factor with separable predual. Denote by $\autM$ the group of automorphisms of M endowed with the $u$-topology, i.e., a sequence of automorphisms $\alpha _{n}$ converges to $\alpha$\, if and only if\, $\|\varphi\circ\alpha _{n}-\varphi\circ\alpha\|\rightarrow 0$ for all $\varphi\in M_{*}$. 

The definition of the Connes invariant $\chi (M)$ involves three normal subgroups of the group of automorphisms $\autM$. For a unitary $u$ in $M$ denote by $\ad _M (u)$ the inner automorphism of $M$ defined by  $\ad _M (u)(x)=uxu^{*}$, for any $x$ in $M$. Let $\innerM$ be the subgroup of $\autM$ formed by all inner automorphisms. Then $\innerM$ is normal in $\autM$ since $\alpha\ad _M (u)\alpha ^{-1}=\ad _M\alpha(u)$ for every $\alpha\in\autM$. Let $\chiusoM$ denote the closure of the group $\innerM$ in the $u$-topology. The group $\innerM$ of inner automorphisms and the group $\chiusoM$ of approximately inner automorphisms are two of the groups involved in the definition of the Connes invariant. 

We restrict now our attention to II$_{1}$ factors. Recall that if 
$\tau$ denotes the unique faithful trace on $M$, then $M$ inherits an $L^{2}$-norm from the inclusion $M\subset L^{2}(M)$, given by $\|x\|_{2}=\tau(x^{*}x)^{1/2}$ for all $x\in M$. Note also that for a II$_{1}$ factor $M$, the 
$u$-topology on $\autM$ is equivalent to the topology for which a sequence of automorphisms $\alpha _{n}$ on $M$ converges to $\alpha$ iff $\displaystyle\lim _{n\rightarrow\infty} \|\alpha _{n}(x)-\alpha (x)\|_{2}\rightarrow 0$. Since our study of automorphisms is simplified in the II$_1$ case, we will always assume that our factors are II$_{1}$, unless otherwise specified. 

The last group of automorphisms involved in the definition of the Connes invariant is formed by the centrally trivial automorphisms of M.

 Given $a,b\in M$ set $[a,b]=ab-ba$. We say that a bounded sequence $(x_{n})_{n\geq 0}$ in $M$ is central if  $\displaystyle\lim_{n \rightarrow \infty}{\|[x_n, y]\|_{2}}=0$ for every $y\in M$.               

\begin{definition}
An automorphism $\alpha \in \autM$ is centrally trivial if for any central sequence $(x_n)$ in $M$ we have $\displaystyle\lim_{n \rightarrow \infty}{\| \alpha (x_n)-x_n \|_{2}}=0.$
\end{definition}

We denote by $\CtM$ the set of centrally trivial automorphisms of $M$, which is a normal subgroup of $\autM$.

Let $\outM=\displaystyle\frac{\autM}{\innerM}$ be the group of outer automorphisms of $M$, and denote by $\xi :\autM\to\outM$ the quotient map. The Connes invariant was introduced by Connes in {\cite{Connes6} (see also \cite{Connes4}).

\begin{definition}
Let $M$ be a II$_{1}$ factor with separable predual. The Connes invariant of $M$ is the abelian group
$$
\chi(M) =\frac{\CtM\cap\chiusoM}{\innerM}\subset\outM.
$$
\end{definition}

Note that the hyperfinite II$_{1}$ factor $R$ has trivial Connes invariant since $\mbox{Ct (R)}=\innerR$.

Next we want to define a particular class of central sequences, the hypercentral sequences. 

\begin{definition}
A central sequence $(x_{n})$ is hypercentral if 
$\displaystyle \lim_{n\rightarrow\infty}\|[x_n , y_n]\|_{2}=0$ for every 
central sequence $(y_{n})$ in $M$. 
\end{definition}

Let $\omega$ be a free ultrafilter over $\mathbb{N}$ and $M$ a II$_{1}$ factor. Denote by  $\ell ^{\infty}(\mathbb{N}, M)$ the algebra of bounded sequences in $M$, and by $C_{\omega}$ the subalgebra of the bounded sequence $(x_{n})_{n\in\mathbb{N}}$ in $M$ with $\displaystyle\lim_{n\rightarrow\omega}\|[x_{n},y]\|_2 =0$, for all $y\in M$. Let $\mathfrak{I}_{\omega}$ be the subalgebra of $\ell ^{\infty}(\mathbb{N}, M)$ consisting of the  sequences for which $\displaystyle\lim _{n\rightarrow\omega}\|x_{n}\|_{2}=0$. Set 
\[M^{\omega}=\frac{\ell ^{\infty}(\mathbb{N}, M)}{\mathfrak{I}_{\omega}}\;\text{ and }\; M_{\omega}=\displaystyle\frac{C_{\omega}}{\mathfrak{I}_{\omega}\cap C_{\omega}}\]
Then, $M^{\omega}$ and $M_{\omega}$ are finite von Neumann algebras and 
\begin{equation}
\label{purple}
M_{\omega}=M^{\omega}\cap M^{\prime}.
\end{equation}

\begin{remark} 
\label{uno}
   By \cite{McDuff} the existence of non-trivial hypercentral sequences is 
   equivalent to the non-triviality of the center of $M_{\omega}$ for some 
   (and then for all) free ultrafilter $\omega$. 
\end{remark}

\begin{definition}  
    A $\mathbb{Z}_{k}$-kernel on a von Neumann algebra $M$ is an automorphism 
    $\alpha\in\autM$ such that there exists a unitary $U$ in $M$ with 
    the property $\alpha ^{k} =\ad _M U$. 
\end{definition}

Note that if $\alpha$ is a $\mathbb{Z}_{k}$-kernel then $\alpha (U)=\lambda U$, 
for $\lambda$ a k-th root of unity. If $\lambda\neq 1$ we say that 
$\alpha$ has obstruction $\lambda$ to lifting, meaning that the homomorphism $\varphi :\mathbb{Z}_{k}\longrightarrow \outM$ defined by $\varphi (1)=[\alpha]$ cannot  be lifted to an homomorphism $\Phi :\mathbb{Z}_{k}\longrightarrow\autM$ with $\Phi (1)=\alpha$.

We conclude this section by defining an invariant $\varkappa (\theta )$ for 
any element  $\theta$ in $\chi (M)$.  

\begin{definition} 
   Let $M$ be a II$_1$ factor without non-trivial hypercentral sequences 
   and take $\theta\in\chi (M)$. Let $\phi$ be an automorphism in $\CtM\cap\chiusoM$ with 
   $\xi (\phi )=\theta$, and $u_{n}$ a sequence of unitaries such that $\phi=\displaystyle\lim _{n\rightarrow\infty}\Ad u_{n}$. Since the sequence $(u_{n}^{*}\phi (u_{n}))_{n\geq 0}$ is hypercentral, there exists a sequence of scalars $(\lambda _{n})_{n\geq 0}$ with the properties that $\displaystyle\lim_{n\rightarrow\infty}\| \phi (u_n )-\lambda _{n}u_{n}\|_{2}=0$, and $(\lambda _{n})_{n\geq 0}$ converges to some $\lambda _{\phi}\in\mathbb{T}$ (Lemmas 2.1 and 2.2 in \cite{Jones1}). Hence,
   \[\lim_{n\rightarrow\infty}{\| \phi (u_n )-\lambda_{\phi }u_n\|}_{2}=0 \] 
and the Jones invariant is $\varkappa (\theta)=\lambda_{\phi }$. 
\end{definition}

Jones proved in \cite{Jones1} that this definition makes sense 
(i.e. $\varkappa (\theta )$ does not depend on the choice of $\phi$ or $u_n$) 
and that $\varkappa$ is a conjugacy invariant, meaning that if $\alpha,\beta$ belong to $\CtM\cap\chiusoM$ and there exists $\psi\in\autM$ such that 
$\psi\alpha\psi ^{-1}=\beta$, then $\varkappa(\xi (\alpha))=\varkappa (\xi 
(\beta))$.

\section{Preliminaries}

Let $M$ be a factor with separable predual. $M$ is said to be full if $\innerM$ is closed in $\autM$ with respect to the $u$-topology. Obviously
all type $I$ factors are full, while the hyperfinite factor $R$ provides an 
example of a II$_{1}$ factor which is not full since $\chiusoR=\autR\neq
\innerR$. 

\begin{remark} 
\label{due}
   For an arbitrary factor, being full is equivalent to having no non-trivial
   central sequence (see \cite{Connes2}). 

\addvspace{\medskipamount}
\end{remark}
The following result, due to Connes, is an easy consequence of Lemma 4.3.3 in \cite{Sakai} and Corollary 3.6 in \cite{Connes1}. Some of the arguments used in the proof can be found in \cite{Jones2}.

\begin{lemma} 
\label{second}
  Let $G$ be a discrete group containing a non-abelian free group and let 
  $\tau$ be the usual trace on $\mathfrak{L}(G)$. Then $\mathfrak{L}(G)$ is full. 
\end{lemma}
 {\bf Proof} Set $E=G-\{e\}$, where $e$  denotes the identity element in $G$.
Let $g_{1},\, g_{2}$ be two generators of the free group and $F=\{g\in E\, 
  |\, g=g_{1}\tilde{g},\;\tilde{g}\in E\}$. Take $x\in \mathfrak{L}(G)$. Then 
  $x$ can be expressed as $x=\displaystyle\sum_{g\in G}{\lambda_g \delta_g}$ and the 
  function $f:G\to \mathbb{C}$ defined by $f(g)=\lambda_g$ belongs to $l^2(G)$.

  For such $f$ we have that 
  \[ \sum_{g\in E}{|f(g)|^2}=\|x-\tau (x)\|_2 ^2\,\text{ and }\, 
  \sum_{g\in G}{|f(g_i g g_i ^{-1})-f(g)|^2} =\|[x,\delta_{g_i}
  ]\|_{2}^{2}.\] 
  Now if $(x_n )$ is a central sequence in $\mathfrak{L}(G)$ then 
  $\| [x_{n},\delta_{g_i}] \|\rightarrow 0$ for $n\rightarrow\infty$, 
  so we can apply Lemma 4.3.3 in \cite{Sakai} and conclude that 
\begin{equation}
\label{puro}
\lim _{n\rightarrow\infty}\|x_n -\tau (x_n )\|_2 =0
\end{equation}

Let $\alpha$ be any automorphism in $\overline{\operatorname {Int}(\mathfrak{L}(G))}$ and choose a sequence of unitaries $(u_{n})_{n}$ such that $\alpha=\displaystyle\lim _{n\longrightarrow\infty} \Ad (u_{n})$. Since $(u_{n}^{*}u_{n+1})_{n\geq 0}$ is a central sequence in $\mathfrak{L}(G)$, by (\ref{puro}) there exists $\lambda _{n}\in\mathbb{T}$ such that 
$$
\|u_{n}^{*}u_{n+1}-\lambda _{n}1\|_{2}<\frac{1}{2^n}.
$$
Set $\displaystyle v_{n}=\left (\prod _{i=1}^{n}{\lambda _{n}}\right )u_{n+1}$. Then $(v_{n})_{n\in\mathbb{N}}$ is a Cauchy sequence so it converges to some $t$ in $\mathfrak{L}(G)$. Since
$$
\Ad (t)=\lim _{n\rightarrow\infty}\Ad (v_{n})=\lim _{n\rightarrow\infty}\Ad (u_{n+1})=\alpha
$$
we have that $\alpha\in\operatorname {Int}(\mathfrak{L}(G))$.
\qvd

Using the following remark we can conclude that not only the free group factors are full, but also the interpolated ones.

\begin{remark}
If $N\subseteq M$ is an inclusion of II$_{1}$ factors and $p$ is a projection in $N$, then $p(N^{\prime}\cap M)p=pN^{\prime}p\cap pMp$. 
\end{remark}

One inclusion of the previous remark is obvious. The other one is proved using the following argument due to S. Popa. Take any element $z\in pN^{\prime}p\cap pMp$ so that $z=x'p$ with $x^{\prime}\in N^{\prime}$. Let $q$ a maximal projection in $N$ with the properties that $p\leq q\leq 1$ and $x^{\prime}q\in M$. To show that $q=1$, suppose that $1-q\neq 0$. Then $(1-q)Np\neq 0$, so using the polar decomposition of a non-zero element in $(1-q)Np$ we can find $0\neq v\in N$ such that $v^{*}v\leq p$ and $vv^{*}\leq 1-q$. Thus $x^{\prime}vv^{*}=vx^{\prime}v^{*}=vpx^{\prime}pv^{*}\in M$ and $x^{\prime}(q+vv^{*})\in M$, contradicting the maximality of $q$.

Taking $N=A$ and $M=A^{\omega}$ in the previous remark, and using (\ref{purple}) we obtain that the compression of a full factor is also full.

\begin{remark}
\label{remark3.4}
Let $A$ be a II$_{1}$ factor with separable predual and $p$ a projection in $A$. Then $A$ is full if and only if $pAp$ is full. In particular, any interpolated free group factor $\inff$, with $t>1$, is full. 
\end{remark}

\begin{proposition}  
\label{sixth}
   Let $\inff$, for $t\in\mathbb{R}$ and $t>1$, be any interpolated free group factor, and denote by $R$ the hyperfinite II$_{1}$ factor. Then $\mathfrak{L}(F_{t})\otimes R$ has no non-trivial hypercentral sequences.
\end{proposition}
{\bf Proof}  
   We start by proving the result for $\mathfrak{L}(G)\otimes R$, where $G$ is a discrete group containing a non-abelian free group. First we want to 
   show that any central sequence in $\mathfrak{L}(G)\otimes R$ has the form 
   $(1\otimes x_{n})_{n\geq 0}$, for a central sequence $(x_{n})_{n\geq 0}$ in $R$.

   Denote by $g_{i}$, for $i=1,2$ the generators of $F_{2}\subseteq G$. By the proof of Lemma \ref{second}, $\mathfrak{L}(G)$ satisfies the hypothesis 
   of Lemma 2.11 in \cite{Connes2} with $Q_1=\mathfrak{L}(G)$, $Q_2=R$ and $b_{i}
   =\delta _{g_{i}}$. Therefore, we can 
   apply the above mentioned lemma to any central sequence $(X_n)_{n\geq 0}$ in  $\mathfrak{L}(G)\otimes R$ to obtain that $\displaystyle\lim_{n\rightarrow\infty}\|X_n-(\tau\otimes 1)(X_n)\|_2=0$. 
   Since $(\tau\otimes 1)(X_n)\in\mathbb{C}\otimes R$, this implies that 
   $X_n=1\otimes x_n$ for a central sequence $(x_{n})_{n\geq 0}$ in $R$. 

   Now suppose $(Y_{n})_{n\geq 0}$ is a hypercentral sequence in $\mathfrak{L}(G)\otimes R$. Since 
   $(Y_{n})$ is central it has the form $Y_n=1\otimes y_n$, for a hypercentral 
   sequence $(y_{n})$ in $R$. So we need only to prove that $R$ has 
   no non-trivial hypercentral sequences.

   This follows immediately from Remark \ref{uno} and Theorem 15.15 in \cite{Kawi}.

   In the case of the factor $\pippo$, we can find  an integer $k>1$ and a projection $p$ in $\mathfrak{L}(F_{k})\otimes R$ such that $\pippo\cong p(\mathfrak{L}(\mathbf{F} _{k})\otimes R)p$. Obviously $p$ belongs to $(\mathfrak{L}(\mathbf{F} _{k})\otimes R)' _{\omega}$. Therefore, by Remark \ref{uno} it is enough to show that given a II$_{1}$ factor $M$ and a projection $p\in M_{\omega}^{\prime}$, $(pMp)_{\omega}$ is a factor if and only if $M_{\omega}$ is a factor. This is an immediate consequence of the equality $(pMp)_{\omega}=pM_{\omega}p$.

\qvd  

\begin{lemma}  
\label{ninth}
   If $\alpha\in\Ctpi$ then $\alpha=\Ad z (\nu\otimes id)$, for some unitary 
   $z\in\pippo$ and an automorphism $\nu$ of $\inff$.

\end{lemma}
   {\bf Proof}
   Let $(K_n)_{n\in\mathbb{N}}$ be an increasing sequence of finite 
   dimensional subfactors of $R$ generating $R$, and denote by $R_n=K_n^{'}\cap R$ the 
   relative commutant of $K_n$ in $R$.

   Set $L_n=1\otimes R_n\subset\pippo$. Then there exists an $n_0$ such that 
   for all $x\in L_{n_0}$,\ \ $\|x\| \leq 1$ one has $\| \alpha(x)-x\|_2\leq
   \frac{1}{2}$. In fact, otherwise it would exist a uniform bounded sequence 
   $(x_{n})$, $x_{n}\in L_{n}$ and $\|x_n\|\leq 1$, such that $\| \alpha(x_n)-
   x_n\|_2 >\frac{1}{2}$. But $(x_{n})$ is a central sequence in $\pippo$ 
   because for each $m$ and $n\geq m$, $x_n$ commutes with $\inff\otimes
   K_m$, so we get a contradiction. 

   By Lemma 3.3 in \cite{Connes1}, up to inner automorphism, $\alpha$ is of the 
   form $\alpha _1\otimes 1_{R_{n_0}}$ where $\alpha _1$ is an automorphism 
   of $\inff\otimes K_{n_0}$. Set $F=1\otimes K_{n_0}$. Then F is a type
   $I$\ subfactor of $\inff\otimes K_{n_0}$.
 
   Applying [Lemma 3.11, \cite{Connes1}] to $\alpha_1 $\ we obtain that 
   $\alpha_1 |_{1\otimes K_{n_0}}=\Ad V |_{1\otimes K_{n_0}}$. This implies 
   that $\alpha=\Ad z(\nu\otimes 1)$ for some automorphism $\nu$ of 
   $\inff$.

\qvd

\section{The subfactor construction of interpolated free group factors}

In this section we use Voiculescu's random matrix model for free 
group algebras \cite{Vocu}, to show  that the crossed product $\crossprod$ can be
realized as the enveloping algebra of an inclusion of interpolated free group factors $A\subset B$. 

For this purpose we first  give an explicit construction of $\crossprod$, 
by providing models for the II$_{1}$ factors $\freefactor$ and $R$. 

Let $\{X_1, X_2, X_3\}$ be a free semicircular family and $u=\displaystyle\sum_{j=1}^{3}{e^{\frac{2\pi ij}{3}} e_j}$ a unitary whose spectral projections $\{e_{j}\}_{j=1} ^{3}$ have trace 
$\frac{1}{3}$. Assume also that $\{u\}^{\prime\prime}$ is free with respect to 
$\{X_1, X_2, X_3\}^{\prime\prime}$. Then $\freefactor$ can be thought as the von Neumann 
algebra generated by $\{X_1, X_2, X_3, u\}$ as in \cite{Vocu} and 
\cite{Radulescu1}.

The model for $R$ is outlined in the following lemma. It is analogous to the 
construction given in the case of $\mathbb{Z}_{2}$ by R\v{a}dulescu \cite{Radulescu2}. We include it here for the sake of completeness.

\begin{lemma}
\label{tenth}
  Given a von Neumann algebra $M$, let $(U_k)_{k\in\mathbb{Z}}$ be a family of unitaries in M of order $9$. Assume that each $U_{k}$ has a decomposition of the form $\displaystyle U_{k}=\sum _{j=1} ^{9}\jonesinv U_{k}^{(j)}$, where each spectral projection $U_{k}^{(j)}$ has trace $\frac{1}{9}$. Let $g=\displaystyle\sum_{j=1}^{3}{e^{\frac{2\pi i j}{3}}g_j}$ be a unitary in $M$ of order $3$ whose spectral projections $\{g_j\}_{j=1}^{3}$ have trace $\frac{1}{3}$. Suppose that the following relations hold between the $U_{k}$'s and $g$:
\begin{enumerate}
 \item [(i)]
    $U_k gU_k^* =\obstruconj\, g\;\text{ if }k=0,-1,\;
    \textstyle{ while }\; U_k gU_k ^* =g\;\text{ if }\;k\in\mathbb{Z}\backslash\{0,-1\}$
 \item [(ii)]
   $U_{k} U_{k+1}U_{k}^{*}=\jonesinv\, U_{k+1},\quad\text{for } k\in\mathbb{Z}$, \item [(iii)]
   $U_{i} U_{j}=U_{j} U_{i}\quad\text{ if}\quad |i-j|\geq 2$.
\end{enumerate}   
The algebra generated by the $U_{k}$'s and $g$ is endowed with a trace defined by $\tau (m)=0$, for each non-trivial monomial $m$ in these unitaries.
Set
$$
R_{-1}=\{gU_{0}^{3}, U_{1}, U_{2},\hdots\}^{\prime\prime}\subset\{g,U_{0}, U_{1}, U_{2},\hdots\}^{\prime\prime}=R_{0}.
$$
This defines an inclusion of type II$_{1}$ factors of index $9$, such that $R_{-1}'\cap R_{0}=\{g\}^{\prime\prime}$.

Let $\theta=\ad _{R_{-1}}(U_0)$. Then $\theta$ is a outer automorphism of
$R_{-1}$ of order $9$ with outer invariant $(3,\obstruction )$. Moreover, 
$R_0$ is equal to the crossed product $\basic$.

Also, the Jones tower for the inclusion $R_{-1}\subset R_0$ is given by
$$
R_{-1}\subset R_{0}\subset R_{1}\subset\cdots\subset R_{k-1}\subset R_{k}\subset R_{k+1}\subset\cdots ,
$$
where $R_k =\{gU_{-1}^3\cdots U_{-k}^{3}, U_{-k}, U_{-k+1},\hdots\}^{\prime\prime}$ for $k\geq 1$. 
\end{lemma}

{\bf Proof}
   The properties of the family $(U_k)_{k\in\mathbb{Z}}$ and of the unitary
   $g$ imply immediately that $U_{0}xU_{0}^{*}\in R_{-1}$ for every $x\in R_{-1}$. Therefore, $\theta=Ad (U_{0})$ defines an automorphism of $R_{-1}$. Since $g$ commutes with $R_{-1}$, $\theta ^{3}=Ad_{R_{-1}}(U_0 ^3)=Ad_{R_{-1}}(gU_0 ^3)$ belongs to $\operatorname{Int}(R_{-1})$. Moreover $\theta (gU_{0}^{3})=\obstruconj gU_{0}^{3}$, so $\theta$ has outer invariant $(3,\obstruconj )$.

Obviously, any monomial in $R_{0}$ can be written using only one occurrence of $U_{0}$ to some power because of the relations between the generators of $R_{0}$. Moreover, by definition, the trace on the algebra generated by the $\{U_{k}\}_{k\in\mathbb{Z}}$ and $g$ (and therefore on its subalgebra $R_{0}$) is compatible with the usual trace defined on the crossed product $\basic$, so that $R_{0}=(R_{-1}\cup\{U_{0}\})^{\prime\prime}=\basic$. 

To show that $R_{-1}^{\prime}\cap R_{0}\subset \{g\}^{\prime\prime}$, write any element $x\in R_{-1}^{\prime}\cap R_{0}$ as $\displaystyle x=\sum _{k=0}^{8}a_{k}U_{0}^{k}$. It is easy to check that $x$ belongs to $R_{-1}^{\prime}$ if and only if $a_{0}\in\mathbb{C}$, $a_{3}$ and $a_{6}$ are multiples of $g^{2}U_{0}^{6}$ and  $gU_{0}^{3}$, respectively, and all the other $a_{k}$'s are zero. The other inclusion follows immediately from the relations verified by the $U_{k}$'s and $g$, thus $R_{-1}^{\prime}\cap R_{0}= \{g\}^{\prime\prime}$.

Note that $\ad _{R_0}(U_{-1})(U_0)=\jonesinv U_0$, while $\ad _{R_0}(U_{-1})(x)=x$ for all $x\in R_{-1}$. Hence, $\ad _{R_0}(U_{-1})$ implements the dual action of $\luna$ on the crossed product $\basic$, and the next step in the Jones tower for the inclusion $R_{-1}\subset R_0$ is given by 
$$
R_1 =\{U_{-1},g,U_{0},U_{1},\hdots\}^{\prime\prime}.
$$
Similarly, the other steps in the Jones constructions are obtained by adding the unitaries $U_{-2},\, U_{-3}, \hdots$, so that the k-th step is given by
\[R_{k} =\{U_{-k},U_{-k+1},\hdots ,U_{-1}, g, U_{0}, U_{1},\hdots\}^{\prime\prime}.\]
\qvd

Observe that to construct unitaries with the properties in the statement it is enough to consider the Jones tower for an inclusion of the form $R\subset R\rtimes_{\beta}\mathbb{Z}_9$, where $\beta$ is an automorphism of order $9$ with outer invariant $(3,\obstruconj )$. We choose the unitary $g$ of order $3$ between the elements of the relative commutant.  The unitaries implementing the crossed product in the successive steps of the basic construction will satisfy the desired relations.  

The next step is to give a concrete realization of the crossed product \linebreak $\M=\crossprod$, with $\gamma$ a $\gruppo$--action on $\tenpro$.

Using the model $\freefactor=\{X_1, X_2, X_3, u\}^{\prime\prime}$, where the $X_{i}$'s are semicircular elements, $u$ is a unitary of order $3$ and $\{X_{i},u \mid \, i=1\hdots 3\}$ is a free family, we define the automorphism $\alpha$ on $\freefactor$ by: 
\begin{itemize}
     \item []
          $\alpha (X_i)=X_{i+1}$, for $i=1,2$
     \item[]$\alpha(X_3)=uX_1 u^*$, and
     \item []
          $\alpha (u)=\obstruction u$.
\end{itemize}
Since $\alpha ^{3}=\Ad u$, $\alpha$ is a $\gruppo$--kernel with obstruction 
$\obstruction$ to lifting.

For the automorphism $\beta$ on the hyperfinite II$_{1}$ factor we use the model of Lemma \ref{tenth}: $R\cong R_0=\{g, U_0, U_1,\hdots\}^{\prime\prime}$ and $\beta =\ad _{R_0}(U_{-1})$, with $\beta ^3 =\ad _{R_0}(g)$ and 
$\beta (g)=\obstruconj g$. 
   
Observe that $\alpha\otimes\beta\in\operatorname{Aut}\pten$ has outer period $3$ and obstruction to lifting $1$, so it can be perturbed by an inner automorphism to obtain a $\gruppo$-action on $\tenpro$. This action is defined by
\begin{displaymath}
\gamma=\left (\ad _{\left [\tenpro\right ]} W\right)(\alpha\otimes
\beta ),
\end{displaymath}
where $W$\ is any cube root of $u^{*}\otimes g^{*}$ which is fixed by the automorphism 
$\alpha\otimes\beta$. For example, if $\delta=\jonesinv$ and $\{e_{i}\}_{i=1} ^{3}$, $\{g_{j}\}_{j=1} ^{3}$ denote the spectral projections of $u$ and $g$, respectively, take 
\begin{equation}
\label{proiezione}
W = \delta E_1+\delta ^{2} E_2+E_3,
\end{equation}
where 
\begin{equation}
\label{Wproj}
E_{l}=\displaystyle\sum_{\substack{i,j=1,\hdots , 3,\\ i+j\equiv l\mod 3}}e_{i}\otimes g_{j}, \text{ for }l=1\hdots , 3.
\end{equation} 

Note that $\alpha$ acts on the spectral projections of $u$ as $\alpha (e_{i})=e_{i-1}$ for $i=2,3$ and $\alpha (e_{1})=e_{3}$, while $\beta$ acts on the spectral projections of $g$ as $\beta (g_{j})=g_{j+1}$ for $j=1,2$ and $\beta (g_{3})=g_{1}$. Hence, $\alpha\otimes\beta$ fixes $W$. 
 
\begin{observation} 
\label{observ}   
   Note that $W$ belongs to the center of the fixed point algebra of $\alpha
   \otimes\beta$ since for any element $z$ in the fixed point algebra we 
   have 
   \begin{equation*}
   z\,(u\otimes g)=(\alpha\otimes\beta )^{3}(z\,(u\otimes g))=\Ad (u\otimes g)
   (z\,(u\otimes g))=(u\otimes g)\,z. 
   \end{equation*}
\end{observation}

We can now prove our main theorem, using an argument similar to the one used by R\u{a}dulescu in \cite{Radulescu2}. We show that $\M=\crossprod$ is the enveloping algebra of an inclusion $A\subset B$ of interpolated free group factors. We divide the proof into two parts, first proving that $A$ is isomorphic to the interpolated free group factor $\mathfrak{L}\left (\mathbf{F} _{\frac{35}{27}}\right )$.   

\begin{proposition}
\label{proposition4.3}
   Let $v$ be the unitary implementing the crossed product $\M=\crossprod$. Consider the von Neumann subalgebra 
$$A=\{X_{i}\otimes 1, u\otimes 1, 1\otimes g, v| i=1,\hdots ,3\}^{\prime\prime}\subset \M,$$ 
endowed with a trace with respect to which $\{X_{i}\otimes 1, u\otimes 1| i=1,\hdots ,3\}$ is a free family, $\{X_{i}\otimes 1\}_{i=1}^{3}$ are semicircular elements, and $u\otimes 1, 1\otimes g, v$ are unitaries of order 3 with spectral projections of trace $\frac{1}{3}$. 
 
Moreover, assume that the following relations are satisfied by the elements of $A$ 
   \begin{itemize}
\item[i)] $[1\otimes g,X_{i}\otimes 1]=0$ for $i=1,\hdots , 3$, and $[1\otimes g,u\otimes 1]=0$. 
\item[ii)] $v(u\otimes 1)=\obstruction\, (u\otimes 1)v$
\item[iii)] $v(1\otimes g)=\obstruconj\, (1\otimes g)v$
\item[iv)] $\Ad v(X_{i}\otimes 1)=\Ad W \circ\alpha(X_{i}\otimes 1)$, where $W=\jonesinv E_{1}+e^{\frac{4\pi i}{9}}E_{2}+E_{3}$ as in (\ref{proiezione}).
\end{itemize}
For any monomial $m$ in the variables $\{X_{i}\otimes 1, u\otimes 1, 1\otimes g| i=1,\hdots , 3\}$, suppose that the trace $\tau$ on the algebra $A$ verifies the following properties: 
   \begin{enumerate} 
   \item [(v)] 
   $\tau(mv^{k})=0$ for $k=1,2$,
   \item [(vi)]
   the trace of $m$ in $A$ coincides with its trace as element of the von Neumann algebra $\{X_{i}, u| i=1, \hdots ,3\}^{\prime\prime}\otimes\{g\}^{\prime\prime}$.
   \end{enumerate}
  
   Under these conditions $A$ is isomorphic to $\mathfrak{L}\left (\mathbf{F} _{\frac{35}{27}}\right )$.

\end{proposition}
{\bf Proof} 
   First we realize the algebra $A$ in terms of random matrices \cite{Vocu} and then use Voiculescu's free probability theory to show that $A$ is an interpolated free group factor. The random matrix model we give is a subalgebra of the algebra of $9\times 9$ matrices with entries in a von Neumann algebra.

   Let $D$ be a von Neumann algebra with a finite trace $\widetilde{\tau}$ which contains a family of free elements $\{a_i\}_{i=1}^{18}$, with the property that the elements $\{a_{i}\}_{i=1} ^{9}$ are semicircular, while the others ones are circular. Denote by $(e_{i j})_{i,j=1,\hdots ,9}$ the canonical system of matrix units in 
   $\matrices$. 
   Set $\epsilon=\obstruction$. Using the same notation as before for the spectral projections of $u$ and $g$, we set  
  $$
e_{i}\otimes 1=\displaystyle \sum _{\substack{j=1,\hdots ,9,\\j\equiv i\mod 3}}e_{jj}\in\matrices\subset D\otimes\matrices,\text{ for } i=1,\hdots ,3, 
$$ 
and 
\begin{align*}
1\otimes g_{1}&=e_{11}+e_{55}+e_{99},\\\ 
1\otimes g_{2}&=e_{33}+e_{44}+e_{88}, \\
1\otimes g_{3}&=e_{22}+e_{66}+e_{77}.
 \end{align*}

Therefore, $u\otimes 1=\epsilon\, (e_{1}\otimes1)+\epsilon ^{2}\, (e_{2}\otimes 1)+e_{3}$ and $1\otimes g=\epsilon\, (1\otimes g_{1})+\epsilon ^{2}\, (1\otimes g_{2})+1\otimes g_{3}$ can be written in matrix notation as
%
%
   \begin{equation*}
   u\otimes 1=\begin{pmatrix}
     \epsilon & 0 & 0 & 0 & 0 & 0 & 0 & 0 & 0 \\
     0 & \epsilon ^2 & 0 & 0 & 0 & 0 & 0 & 0 & 0 \\
     0 & 0 & 1 & 0 & 0 & 0 & 0 & 0 & 0 \\
     0 & 0 & 0 & \epsilon & 0 & 0 & 0 & 0 & 0 \\
     0 & 0 & 0 & 0 & \epsilon ^2 & 0 & 0 & 0 & 0 \\
     0 & 0 & 0 & 0 & 0 & 1 & 0 & 0 & 0 \\
     0 & 0 & 0 & 0 & 0 & 0 & \epsilon & 0 & 0 \\
     0 & 0 & 0 & 0 & 0 & 0 & 0 & \epsilon ^2 & 0 \\
     0 & 0 & 0 & 0 & 0 & 0 & 0 & 0 & 1
     \end{pmatrix}
   \end{equation*}
  
  \addvspace{\baselineskip}
  \begin{equation*}
   1\otimes g=\begin{pmatrix}
     \epsilon & 0 & 0 & 0 & 0 & 0 & 0 & 0 & 0 \\
     0 & 1 & 0 & 0 & 0 & 0 & 0 & 0 & 0 \\
     0 & 0 & \epsilon ^2 & 0 & 0 & 0 & 0 & 0 & 0 \\
     0 & 0 & 0 & \epsilon ^2 & 0 & 0 & 0 & 0 & 0 \\
     0 & 0 & 0 & 0 & \epsilon & 0 & 0 & 0 & 0 \\
     0 & 0 & 0 & 0 & 0 & 1 & 0 & 0 & 0 \\
     0 & 0 & 0 & 0 & 0 & 0 & 1 & 0 & 0 \\
     0 & 0 & 0 & 0 & 0 & 0 & 0 & \epsilon ^2 & 0 \\
     0 & 0 & 0 & 0 & 0 & 0 & 0 & 0 & \epsilon
     \end{pmatrix}
   \end{equation*}
   
  \addvspace{\baselineskip}
Moreover, set  
$$ 
v=(e_{1 2}+e_{2 3}+e_{3 1})+(e_{4 5}+e_{5 6}+e_{6 4})+(e_{7 8}+e_{8 9}+e_{9 7}),$$
i.e., 

\begin{equation*}
   v=\begin{pmatrix}
     0 & 1 & 0 & 0 & 0 & 0 & 0 & 0 & 0 \\
     0 & 0 & 1 & 0 & 0 & 0 & 0 & 0 & 0 \\
     1 & 0 & 0 & 0 & 0 & 0 & 0 & 0 & 0 \\
     0 & 0 & 0 & 0 & 1 & 0 & 0 & 0 & 0 \\
     0 & 0 & 0 & 0 & 0 & 1 & 0 & 0 & 0 \\
     0 & 0 & 0 & 1 & 0 & 0 & 0 & 0 & 0 \\
     0 & 0 & 0 & 0 & 0 & 0 & 0 & 1 & 0 \\
     0 & 0 & 0 & 0 & 0 & 0 & 0 & 0 & 1 \\
     0 & 0 & 0 & 0 & 0 & 0 & 1 & 0 & 0
     \end{pmatrix}.
   \end{equation*} 
Note that the unitaries $u\otimes 1$, $1\otimes g$, and $v$ generate a copy of $M_3 (\mathbb{C})\oplus M_3 (\mathbb{C})\oplus M_3 
   (\mathbb{C})\subseteq\matrices$. 
   In addition, with this choice of $u\otimes 1$, $1\otimes g$ and $v$ we obtain that
   \begin{equation*}
   W=\delta ^2Id\oplus Id\oplus\delta\, Id,
   \end{equation*}
   where $\delta=\jonesinv$ and $Id$ denotes the identity of $M_3 (\mathbb{C})$. In matrix 
   notation 
   \begin{equation*}
      W=\begin{pmatrix}
              \delta ^2 & 0 & 0 & 0 & 0 & 0 & 0 & 0 & 0 \\
              0 & \delta ^2 & 0 & 0 & 0 & 0 & 0 & 0 & 0 \\
              0 & 0 & \delta^2 & 0 & 0 & 0 & 0 & 0 & 0 \\
              0 & 0 & 0 & 1 & 0 & 0 & 0 & 0 & 0 \\
              0 & 0 & 0 & 0 & 1 & 0 & 0 & 0 & 0 \\
              0 & 0 & 0 & 0 & 0 & 1 & 0 & 0 & 0 \\
              0 & 0 & 0 & 0 & 0 & 0 & \delta & 0 & 0 \\
              0 & 0 & 0 & 0 & 0 & 0 & 0 & \delta & 0 \\
              0 & 0 & 0 & 0 & 0 & 0 & 0 & 0 & \delta
     \end{pmatrix}
   \end{equation*}
   Furthermore, $u\otimes 1$, $1\otimes g$ and $v$ satisfy the required conditions
   \begin{align*}
   (u\otimes 1)(1\otimes g) &=(1\otimes g)(u\otimes 1),\\
   v(u\otimes 1)&=\epsilon \, (u\otimes 1) v,\\ 
   v(1\otimes g)&=\overline{\epsilon}\, (1\otimes g)v.
   \end{align*}
   Denote by $\Rea(a)=\frac{1}{2}(a+a^* )$, for $a\in D\otimes\matrices$.
   We set:
   \begin{align*}
X_{1}\otimes 1= & a_1\otimes e_{1 1}+a_2\otimes e_{22}+a_3\otimes e_{33}+a_{4}\otimes 
       e_{44}+a_{5}\otimes e_{55}+a_{6}\otimes e_{66} \\
     & +a_{7}\otimes e_{77}+a_{8}\otimes e_{8 8}+a_{9}\otimes e_{99}+2\Rea(a_{10}\otimes e_{15})+2\Rea(a_{11}\otimes e_{19}) \\
     & +2\Rea(a_{12}\otimes e_{26})+2\Rea(a_{13}\otimes e_{27})+2\Rea(a_{14}\otimes e_{3 4})+2\Rea(a_{15}\otimes e_{38}) \\
     & +2\Rea(a_{16}\otimes e_{48})+2\Rea(a_{17}\otimes e_{59})+2\Rea(a_{18}\otimes e_{67}).
   \end{align*}
   To find $X_2$ and $X_3$, use the relations:
   \begin{equation*}
   X_2\otimes 1=\Ad (W^* v)(X_1\otimes 1)\quad\text{ and }\quad X_3\otimes 1=\Ad (W^* v)(X_2\otimes 1)
    \end{equation*}

  Thus, using matrix notation we have: 
   \begin{equation*}
       X_1\otimes 1=\begin{pmatrix} 
                   a_{1} & 0 & 0 & 0 & a_{10} & 0 & 0 & 0 & a_{11}  \\ 
                   0 & a_{2} & 0 & 0 & 0 & a_{12} & a_{13} & 0 & 0 \\ 
                   0 & 0 & a_{3} & a_{14} & 0 & 0 & 0 & a_{15} & 0 \\ 
                   0 & 0 & a_{14}^* & a_{4} & 0 & 0 & 0 & a_{16} & 0 \\ 
                   a_{10}^* & 0 & 0 & 0 & a_{5} & 0 & 0 & 0 & a_{17} \\ 
                   0 & a_{12}^* & 0 & 0 & 0 & a_{6} & a_{18} & 0 & 0 \\ 
                   0 & a_{13}^* & 0 & 0 & 0 & a_{18}^* & a_{7} & 0 & 0 \\ 
                   0 & 0 & a_{15}^* & a_{16}^* & 0 & 0 & 0 & a_{8} & 0 \\ 
                   a_{11}^* & 0 & 0 & 0 & a_{17}^* & 0 & 0 & 0 & a_{9}
              \end{pmatrix}
   \end{equation*}
   
   \begin{equation*}
      X_2\otimes 1=\begin{pmatrix}
    a_{2} & 0 & 0 & 0 & \jconj ^{2} a_{12} & 0 & 0 & 0 & \jconj a_{13} \\
    0 & a_{3} & 0 & 0 & 0 & \jconj ^{2} a_{14} & \jconj a_{15} & 0 & 0 \\
    0 & 0 & a_{1} & \jconj ^{2} a_{10} & 0 & 0 & 0 & \jconj a_{11} & 0 \\
    0 & 0 & \delta ^{2} a_{10}^{*} & a_{5} & 0 & 0 & 0 & \delta a_{17} & 0 \\
    \delta ^{2} a_{12}^{*} & 0 & 0 & 0 & a_{6} & 0 & 0 & 0 & \delta a_{18} \\
    0 & \delta ^{2} a_{14}^{*} & 0 & 0 & 0 & a_{4} & \delta a_{16} & 0 & 0 \\
    0 & \delta a_{15}^{*} & 0 & 0 & 0 & \jconj a_{16}^{*} & a_{8} & 0 & 0 \\
    0 & 0 & \delta a_{11}^{*} & \jconj a_{17}^{*} & 0 & 0 & 0 & a_{9} & 0 \\
    \delta a_{13}^{*} & 0 & 0 & 0 & \jconj a_{18}^{*} & 0 & 0 & 0 & a_{7} 
   \end{pmatrix}
   \end{equation*}

  \begin{equation*} 
      X_3\otimes 1=\begin{pmatrix}
a_{3} & 0 & 0 & 0 & \jconj ^{4} a_{14} & 0 & 0 & 0 & \jconj ^{2} a_{15} \\
0 & a_{1} & 0 & 0 & 0 & \jconj ^{4} a_{10} & \jconj ^{2} a_{11} &  0 & 0 \\
0 & 0 & a_{2} & \jconj ^{4} a_{12} & 0 & 0 & 0 & \jconj ^{2}a_{13} & 0 \\
0 & 0 & \delta ^{4}a_{12}^{*} & a_{6} & 0 & 0 & 0 & \delta ^{2} a_{18} & 0 \\  
\delta ^{4}a_{14}^{*} & 0 & 0 & 0 & a_{4} & 0 & 0 & 0 & \delta ^{2} a_{16} \\
0 & \delta ^{4}a_{10}^{*} & 0 & 0 & 0 & a_{5} & \delta ^{2} a_{17} & 0 & 0 \\
0 & \delta ^{2} a_{11}^{*} & 0 & 0 & 0 & \jconj ^{2} a_{17}^{*} & a_{9} & 0 & 0 \\
0 & 0 & \delta ^{2} a_{13}^{*} & \jconj ^{2} a_{18}^{*} & 0 & 0 & 0 & a_{7} & 0 \\
\delta ^{2} a_{15}^{*} & 0 & 0 & 0 & \jconj ^{2} a_{16}^{*} & 0 & 0 & 0 & a_{8}
\end{pmatrix}
\end{equation*}
  
\addvspace{1.3\baselineskip}
Obviously $X_1\otimes 1 , X_ 2\otimes 1 ,X_3\otimes 1$ commute with $1\otimes g$. Moreover, the assumption that the family $\{a_i\}_{i=1,\hdots ,18}$\ is free implies that $\{X_{i}\otimes 1, u\otimes 1 |i=1,\hdots ,3\}$ is a free family with respect to the unique normalized trace on $D\otimes\matrices$. In 
addition, $\{X_{i}\otimes 1, u\otimes 1|i=1,\hdots ,3\}^{\prime\prime}$ and $\{g\}^{\prime\prime}$ are independent with respect to this trace. 

  To show that the normalized trace of $D\otimes\matrices$ has the properties (v) and (vi) of the statement, let $m$ be any monomial in the variables $\{X_{i}\otimes 1, u\otimes 1, 1\otimes g| i=1,\hdots ,3\}$, and consider the product $v^k m$ with $k=1,2$. Since $m$ commutes with $1\otimes g$, it 
  must  be of the form 
  \begin{equation*}
  m=\begin{pmatrix}
    * & 0 & 0 & 0 & * & 0 & 0 & 0 & * \\
    0 & * & 0 & 0 & 0 & * & * & 0 & 0 \\
    0 & 0 & * & * & 0 & 0 & 0 & * & 0 \\
    0 & 0 & * & * & 0 & 0 & 0 & * & 0 \\
    * & 0 & 0 & 0 & * & 0 & 0 & 0 & * \\
    0 & * & 0 & 0 & 0 & * & * & 0 & 0 \\
    0 & * & 0 & 0 & 0 & * & * & 0 & 0 \\
    0 & 0 & * & * & 0 & 0 & 0 & * & 0 \\
    * & 0 & 0 & 0 & * & 0 & 0 & 0 & *
    \end{pmatrix}
  \end{equation*}  
  
  If we multiply $m$ by $v$ or $v^2$, we obtain a matrix with zero on
  the diagonal and a few non-zero entries outside the diagonal. This implies
  that $v^k m$ has zero trace for $k=1,2$. Thus, we have built a matrix model 
  for $A$ which satisfies the conditions of the statement. 
  
  One can easily check that $A$ is a factor. In order to prove that $A$ is an interpolated free group factor we reduce $A$ by one of the spectral projections of $g$, and show that the new factor we obtain is an interpolated free group factor. Reduce $A$ by $g_3 =1\otimes e_{22}+1\otimes e_{66}+1\otimes e_{77}$, which has trace $\frac{1}{3}$.

Then $g_3 A g_3$ is generated by:
\begin{enumerate}
  \item [(i)] 
         $a_{2}\otimes e_{22}+2\Rea(a_{12}\otimes e_{26})+2\Rea(a_{13}\otimes 
         e_{27})+a_{6}\otimes e_{66}+ a_{7}\otimes e_{77} \\
         +2\Rea(a_{18}\otimes e_{67})$
  \item [(ii)]
     $a_{3}\otimes e_{22}+2\Rea(\jconj ^{2}a_{14}\otimes e_{26})+
    2\Rea(\jconj a_{15}\otimes e_{27})+a_{4}\otimes e_{66}+a_{8}\otimes e_{77} \\
    +2\Rea(\delta a_{16}\otimes e_{67})$ 
 \item  [(iii)]
    $a_{1}\otimes e_{22}+2\Rea(\jconj ^{4}a_{10}\otimes e_{26})+
    2\Rea(\jconj ^{2} a_{11}\otimes e_{27})+a_{5}\otimes e_{66}+a_{9}\otimes e_{77} \\
    +2\Rea(\delta ^{2} a_{17}\otimes e_{67})$
 \item  [(iv)]
    $\epsilon ^{2}\otimes e_{22}+1\otimes e_{66}+\epsilon\otimes e_{77}$,
\end{enumerate}
  Thus, by the Voiculescu's random matrix model \cite{VoDy} $g_3 A g_3$ is isomorphic
  to the free group factor $\mathfrak{L}(F_{3}*\mathbb{Z}_{3})\cong\freefactor$. 
  This implies that $A$ is also an interpolated free group factor, and using the well-known 
  formula for reduced factors (\cite{Dykema1}, \cite{Dykema2} or \cite{Radulescu1}) 
  we get that $A=\mathfrak{L}\left (\mathbf{F} _{\frac{35}{27}}\right )$. 

\qvd

\begin{theorem} 
\label{main}
  Set  
  $$
  A=\{X_{i}\otimes 1,u\otimes 1, 1\otimes g,v\mid i=1\hdots 3\}^{\prime\prime}\subset B=(A\cup \{1\otimes U_{0}\})^{\prime\prime}.
  $$
  in $\M=\crossprod$. Then $A$ is isomorphic to the interpolated free group factor $\mathfrak{L}\left (\mathbb{F}_{\frac{35}{27}}\right )$, and $B$ 
  is the crossed product of $A$ by a $\luna$-action $\theta$ on $A$, with outer invariant is $(3, \obstruction )$. Furthermore, $\M$ is 
  the enveloping algebra in the basic construction for the inclusion $A\subset B$. 
\end{theorem}

{\bf Proof} First we want to describe how $\Ad (1\otimes U_0)$ acts on the subalgebra $A$ of $\M$. 

  Obviously
  \begin{equation*}
  \Ad (1\otimes U_{0})(X_{i}\otimes 1)=X_{i}\otimes 1\quad\text{ for }i=1,\hdots ,3,
  \end{equation*}
  \begin{equation*}
  \Ad (1\otimes U_{0})(u\otimes 1)=u\otimes 1
  \end{equation*}
  and 
  \begin{equation*}
  \Ad (1\otimes U_{0})(1\otimes g)=\obstruconj (1\otimes g).
  \end{equation*}
In addition, if $\{E_{i}\}_{i=1} ^{3}$ are the projections defined in (\ref{Wproj}), then
\begin{equation}
\label{E}
 \Ad (1\otimes U_{0})(E_{i})=E_{i+1}\text{ for } i=1,\hdots , 3,
\end{equation}
 where $i+1$ is taken mod 3. Using (\ref{proiezione}) and (\ref{E}), together with the relation $\gamma (1\otimes U_{0})=\delta \Ad W(1\otimes U_{0})$, for  $\delta=\jonesinv$, we obtain 
  \begin{align*}
  \Ad (1\otimes U_{0})(v) & =(1\otimes U_{0})(v(1\otimes U_{0}^{*})v^* )v=\bar{\delta}\,\Ad (1\otimes U_0)(W)W^{*} v \\
  & =\bar{\delta}\,(E_{1}+\delta E_{2}+\delta ^{2}E_{3})W^{*}v 
  =\bar{\delta ^{2}}\, (E_1 +E_2 +\obstruction E_3) v.
  \end{align*}

  It follows that $\Ad (1\otimes U_0)$ leaves $A$ invariant. Furthermore, 
  $\ad _{A}(1\otimes U_{0})^{3}$ acts identically on $\{X_{i}\otimes 1, u\otimes 1, 1\otimes g| i=1,\hdots , 3\}''$,   
  while $\ad _{A}(1\otimes U_{0})^{3} (v)=\obstruconj\, v$.

  Set $\theta=\ad _{A}(1\otimes U_{0})$. Then, using the fact that 
  $\ad _{A}(1\otimes g^*)$  acts identically on $\{X_{i}\otimes 1, u\otimes 1, 1\otimes g| i=1,\hdots ,3\}^{\prime\prime}$, and $\Ad v(1\otimes g)=\obstruconj (1\otimes g)$, we conclude 
that 
  $$
  \theta ^{3}=\ad _{A}(1\otimes g^*)\quad\text{ and }\quad\theta(1\otimes g^*)=
  \obstruction\, (1\otimes g^*) .
  $$

  Thus $\theta$ is a $\gruppo$-kernel on $A$ with obstruction $\obstruction$
  to lifting. 

  To complete the proof that $B=A\rtimes_{\theta}\luna$ we need to check that
  any monomial $m$ in the variables $\{X_{i}\otimes 1, u\otimes 1, 1\otimes g, 1\otimes U_0 , v| i=1,\hdots ,3\}$ contains 
only one occurrence  of $1\otimes U_0$ to some power. We also need to 
  verify that any monomial containing $1\otimes U_0$ has zero trace. 

  Since the crossed product $\crossprod$ is implemented by the unitary $v$, 
  any monomial in the variables $\{X_{i}\otimes 1, u\otimes 1,1\otimes g, 1\otimes U_0 , v| i=1,\hdots ,3\}$ 
  has the form $v^{k} m$, where $m$ is an element in $\tenpro$ and 
  $k=0,1,2$. Because $1\otimes U_{0}$ commutes with the elements of $\freefactor\otimes 1$ and 
$\Ad (1\otimes U_{0})(1\otimes g)=\obstruconj\, (1\otimes g)$, it follows that $1\otimes U_{0}$ can appear at most once in $m$ with some power. Thus, 
any monomial in $B$ can be written using only one occurrence of 
  $1\otimes U_0$. 

  Furthermore, by the definition of the trace on the crossed product, any monomial of the form $v^k m$, for $k=1,2$, has zero trace 
  in $B$. Therefore, it is enough to compute the trace of any monomial $m$ in $\tenpro$ containing $1\otimes U_{0}$. Because of the definition 
of the trace on $R_{0}$ (Lemma \ref{tenth}) and hence on $\tenpro$, we can conclude immediately that the trace of $m$ is zero. Therefore 
  $B=A\rtimes_{\theta}\luna$.

  In addition, using an argument similar to the one used to build the model 
  for $R$ (Lemma \ref{tenth}), one can easily verify that the unitaries $1\otimes U_1, 
  1\otimes U_2,\hdots$\, implement the consecutive terms in the iterated basic 
  construction of $A\subset A\rtimes_{\theta}\luna$. This implies that 
  $\crossprod$\ is the enveloping algebra of the above inclusion of 
  factors.

\qvd

\section{The Connes invariant of the crossed product $\crossprod$}

The arguments in this section are analogous to the ones used by Jones in 
\cite{Jones1} to compute the Connes invariant for his example of a factor 
which is anti-isomorphic to itself but has no involutory antiautomorphism. The definition of $K$, $K^{\bot}$ and $L$ given below are due to 
Connes (see \cite{Connes6}).

Let $N$\ be a II$_{1}$ factor without non-trivial hypercentral sequences and
$G$ a finite subgroup of $\autN$ such that $G \cap\chiusoN=\{Id\}$. 
Set $K=G\cap\CtN$ and \linebreak
$K^{\bot}=\{f:G\to\mathbb{T}\, | \text{ f is a homomorphism and } f|_{K}\equiv
1\}$.
 
Let $\xi :\autN\to\outN$ be the usual quotient map and Fint be
the subgroup $\fint\subseteq\autN$. Denote by $\closure$ its closure in 
$\autN$ with respect to the pointwise weak topology, and by $G\vee\CtN$ the subgroup of $\autN$ generated by $G\cup\CtN$. Set $L=\xi ((G\vee\CtN)\cap\closure)\subseteq\outN$. 

Connes showed in \cite{Connes6} that there exists an exact sequence
\begin{equation*}
0\rightarrow K^{\bot}\overset{\partial}{\rightarrow}\chi (W^{*}(N,G))
\overset{\Pi}{\rightarrow}  L\rightarrow 0
\end{equation*}
where $M=W^{*}(N,G)$ denotes the crossed product implemented by the action of $G$ on $N$. 
We briefly describe the maps $\partial$ and $\Pi$ in the exact sequence above. Given an element $x$ in $M$, write it as 
$\displaystyle\sum_{g\in G}{a_{g}u_g }$, with $a_{g}\in N$ and $\ad _{M}u_{g}|_{N}=g$. 
For each $\eta :G \rightarrow \mathbb{T}$ in $K^{\bot}$, define the map $\Delta (\eta):M\longrightarrow M$ by
\begin{displaymath}
\Delta (\eta )\left (\sum_{g\in G}{a_g u_g}\right )= \sum_{g\in G}{\eta (g)
a_{g}u_{g}}.
\end{displaymath}
Then $\Delta (\eta )$ belongs to $\CtM\cap \chiusoM$, and $\partial=\xi\circ\Delta$ is the desired map. 

To see how  $\Pi$ acts on $\chi (M)$, for any element $\sigma\in\chi (M)$ choose an automorphism $\alpha\in\CtM\cap\chiusoM$ such that $\xi (\alpha )=\sigma$. 

The hypothesis $G\cap\chiusoN=\{Id\}$, implies that there exists a sequence of unitaries $(u_{n})_{n\geq 0}$ in $N$, which are fixed by $G$, and a unitary $z$ in $M$ such that $\alpha=\displaystyle \Ad z\lim _{n\rightarrow\infty}{\Ad u_n}$ (Corollary 6 and Lemma 2 in \cite{Jones2}, or Lemma 15.42 in \cite{Kawi}). Set $\psi _{\sigma}\Ad (z^*)\alpha|_{N}\in\autN$.  

One can show that $\psi_{\sigma}\in \closure\cap G\vee\CtN$, and that the composition map $\xi\circ\psi _{\sigma}$ does not depend on the choice of $\alpha$, but only on the class $\sigma=\xi (\alpha)$. Therefore, the map \,$\Pi :\chi(M)\to L$\, given by $\Pi(\sigma )=\xi (\psi_{\sigma })$ is well defined. 

To show that $\Pi$ is surjective, let $\mu$ be any element in $L$ and denote by $\alpha_{\mu}\in\closure\cap (G\vee\CtN)$ a representative of $\mu$, i.e. $\xi (\alpha_{\mu})=\mu$.  $\alpha_{\mu}$ commutes with $G$ since it is
the limit of automorphisms with this property. Hence, the map $\beta_{\mu}$ defined by
\begin{equation}
\label{lifting}
 \beta_{\mu}\left (\sum_{g\in G}{a_{g}u_g }\right )=\sum_{g\in G}{\alpha_{\mu}(a_g)u_g}.
\end{equation}
is an automorphism of $M$. In addition, we have that $\beta_{\mu}\in\CtM\cap\chiusoM$ and $\Pi (\xi (\beta_{\mu}))=\mu$. 

\addvspace{\baselineskip}
\begin{remark}[Jones] 
\label{tre}
   Note that if $u_{n}$ is a sequence of unitaries left invariant by $G$, 
   with the property $\alpha_{\mu}=\displaystyle\lim_{n\rightarrow\infty}{\Ad u_{n}}$ in $\autN$, and  $\beta_{\mu}$ is the map defined in (\ref{lifting}), then 
   $\beta_{\mu}=\displaystyle\lim_{n\rightarrow\infty} {\Ad u_{n}}$ in 
   $\autM$. Hence,
   \begin{equation*}
   \varkappa (\varepsilon (\beta_{\mu}))=\lim_{n\rightarrow\infty}u_{n}^{*}
   \beta_{\mu}(u_{n})=\lim_{n\rightarrow\infty}{u_{n}^{*}\alpha_{\mu}
   (u_{n})}.  
   \end{equation*}  

\addvspace{\baselineskip} 
\end{remark}

Our next goal is to show that if $N=\tenpro$ and $G$ is the subgroup of $\autpten$ generated by $\gamma =\Ad W (\alpha\otimes\beta),$ then $\chi (\M)\cong\luna$, as it was observed by Connes in \cite{Connes4}. 

For all the rest of this paper we will identify $G=\langle\gamma\rangle$ with 
$\gruppo$. Note that by Proposition \ref{sixth}, $\tenpro$ has no non-trivial 
hypercentral sequences, so in order to use the exactness of the Connes
sequence described above, we only need to show that $\gruppo\cap\chiusopten =
\{Id \}$. Obviously it is enough to check that $\gamma\not\in\chiusopten$. This is equivalent to show that $\alpha\otimes\beta\not\in\chiusopten$. 

By \cite[Corollary 3.3]{Connes3} if $\alpha\otimes\beta\in\chiusopten$ then $\alpha\in\chiusofac$ and $\beta\in\chiusoR$. But $\freefactor$ is 
full (Remark \ref{remark3.4}) and $\alpha\not\in\innerfac$, thus we can conclude that $\gruppo\cap\chiusopten =\{Id\}$.

Hence, for the factor $\M=\crossprod$ the sequence 
$$
0 \rightarrow K^{\bot}\overset{\partial}{\rightarrow}\chi (\M)\overset{\Pi}
{\rightarrow} L\rightarrow 0
$$ 
is exact (\cite{Connes6}). In order to compute $\chi (\M)$ we first show that 
$$
K^{\bot}=\gruppo\quad\text{and}\quad L=\gruppo.
$$ 

\begin{lemma}
   The group $K=\gruppo\cap\Ctpten$ is trivial.
\end{lemma}
{\bf Proof} 
   Obviously it is enough to show that the automorphism $\gamma=Ad\, W 
   (\alpha\otimes\beta )$ is not in 
   $\Ctpten$. We have already shown in the proof of Lemma \ref{sixth} 
   that any central sequence in $\tenpro$ has the form 
   $(1\otimes x_n )_{n\in\mathbb{N}}$ for a central sequence $(x_{n}
   )_{n\in\mathbb{N}}$ in $R_{0}$. It follows that $\gamma\in\Ctpten$ if and only 
   if $\beta\in\CtR$. Since for $\epsilon=\obstruction$, $\beta$ is outer 
   conjugate to the automorphism $s_{3}^{\bar{\epsilon}}$ described by Connes   
   in \cite{Connes5}, and $s_{3}^{\bar{\epsilon}}$ does not belong to $\CtR$ \cite[Proposition 1.6]{Connes5}, we conclude that $\beta\not\in\CtR$. 

\qvd
  
\begin{lemma}
\label{eleventh}
    The group $L$ is isomorphic to $\gruppo$, and a generator of $L$ is given by 
    \[ \mu=\xi\left (\gener\right ), \]
    where $U_0$\ is the unitary in $R_{0}$ defined in Lemma \ref{tenth}, such that $\beta (U_{0})=\jonesinv U_{0}$. 

\end{lemma}    
{\bf Proof}
   We need to show that the automorphism $\gener$ belongs to \linebreak $\left (\gruppo\vee\Ctpten\right )\cap\closure$. 
   To prove that it is in $\gruppo\vee\Ctpten$, multiply $\gener$ by 
   $\gamma^{-1}=(\Ad W (\alpha\otimes\beta))^{-1}$ to obtain the automorphism 
   \[ \Ad\left (W^* (1\otimes U_{0}^{*} )\right )(\alpha^{-1}\otimes 
   Id), \]
   which is in $\Ctpten$, since the central sequences in $\tenpro$ have the form $1\otimes x_{n}$, with $x_{n}$ central in $R_{0}$. Thus 
   $\gener\in\gruppo\vee\Ctpten$.

   Next we want to show that $\gener\in\closure$. Therefore, we need to exhibit a 
   sequence $(\tilde{u} _{n})$ of unitaries in $\tenpro$, that are invariant with respect to 
   $\gamma$, and satisfy $\gener =\displaystyle\lim_{n\rightarrow\infty}{\Ad\tilde{u} _n}$.

   Observe that the sequence $(x_n)_{n\in\mathbb{N}}$ of unitaries in $R_{0}$ given by 
   \begin{equation*}
   x_n=\begin{cases}
        U_{0}U_{1}^{*}U_{2}U_{3}^{*}\hdots U_{n}^{*},&\text{ if $n$ is odd}, \\
        U_{0}U_{1}^{*}U_{2}U_{3}^{*}\hdots U_{n},& \text{ if $n$ is even} \\
       \end{cases} 
   \end{equation*}
   has the properties 
   \begin{equation*}
   \beta =\lim_{n\rightarrow\infty}{\Ad x_n}\quad\text{ and }\quad\beta 
   (x_n)=\jonesinv\, x_n.
   \end{equation*}

   Define $u_{n}=U_{0}^{*} x_n$. Obviously 
   \begin{equation*}
   Id\otimes (\Ad U_{0}^{*}\beta) =\lim_{n\rightarrow\infty}\Ad (1\otimes u_n)\quad\text{ and
    }\quad\beta (u_n)=u_n
   \end{equation*} 
   so that 
   \begin{equation}
   \label{quarta} 
   (\alpha\otimes\beta)(1\otimes u_n)=1\otimes u_n. 
   \end{equation}
   In addition, from Observation \ref{observ} it follows that  
   \[ \gamma(1\otimes u_{n})=\Ad W (\alpha\otimes\beta)(1\otimes
   u_{n})=1\otimes u_{n}. \]
   Thus, 
  \[\gener\in\left (\gruppo\vee\Ctpten\right )\cap\closure.\]

   Lastly we want to prove that the order 3 element $\mu =\xi\left (\gener\right )$ generates 
   $L$. Thus we need to show that any element $\varphi$ in $\left (\gruppo\vee\Ctpten\right )\cap\closure$, is of the form $\mu ^{n}\Ad w$, for some unitary $w$ in $\tenpro$ and $n=0,\hdots ,2$. Since $\varphi\in\gruppo\vee\Ctpten$, there exists $k\in\{0,1,2\}$ such that $\varphi\gamma ^{k}$ is 
   centrally trivial. By Proposition 3.6, it follows that there exists a unitary  $z\in\tenpro$ and an automorphism $\nu\in\autfac$ such that 
   $\varphi\gamma ^{n}=\Ad z(\nu\otimes id)$. Therefore, 
   \[ \varphi=\Ad x(\nu\alpha^{-n}\otimes\beta^{-n}), \]
   with $x=z(\nu\otimes Id)(W^{n})^{*}\in\tenpro$.

   Since $\varphi\in\closure\subseteq\chiusopten$ and $\freefactor$ 
   is full (Remark \ref{remark3.4}), by \cite[Corollary 3.3]{Connes3} there exists a unitary $w$ in $\freefactor$ such 
   that $\nu\alpha^{-n}=\Ad w$. This implies that 
   $\varphi=\Ad x' (Id\otimes\beta ^{-n})$, where $x'=x\,(w\otimes 1)$. Thus 
   $\varphi$ differs from a power of $\gener$ only by an inner automorphism. Hence 
   $\gener$ generates $L$.

   \qvd

\addvspace{\baselineskip}
Note that if $u_{n}$ is the sequence of unitaries defined in the previous lemma 
\[ (\gener)(1\otimes u_{n})=1\otimes U_{0}^{*}u_{n}U_0. \]
But $\displaystyle\lim_{n\rightarrow\infty}{\Ad u_n }=\Ad U_{0}^*\,\beta$ and 
$\beta(U_{0})=\jonesinv\, U_{0}$, so 
\[ \lim_{n\rightarrow\infty}{u_{n}^{*}U_{0}^{*}u_{n}}=\Ad U_{0}(\beta^{-1}(U_{0}^{*}))=
\jonesinv\, U_{0}^{*} \] 
and
\begin{equation*}
\lim_{n\rightarrow\infty}{(1\otimes u_{n}^{*})(\gener )(1\otimes u_n )}=\jonesinv.
\end{equation*}

Thus, in view of Remark \ref{tre} we obtain the following:

\begin{remark}
\label{add}   
   Let $\mu=\xi\left (\gener\right )$ be the generator of $L$ as in the previous lemma. If $\beta_{\mu}\in\CtM\cap\chiusoM$ is the automorphism described in
   equation (\ref{lifting}), and $\sigma=\varepsilon (\beta_{\mu})$, then 
   $\varkappa (\sigma)=\jonesinv$.

\addvspace{\medskipamount}
\end{remark}
\begin{theorem}
\label{twelve}  
   Let $\M$ denote the crossed product $\crossprod$. Then $\chi (\M)=\luna$.

\end{theorem}
{\em Proof}. 
  By Connes \cite{Connes6} the sequence
  \begin{displaymath}
  0 \rightarrow\gruppo\rightarrow\chi (\M)\rightarrow\gruppo\rightarrow 0 
  \end{displaymath}
  is exact. Moreover, according to (\ref{lifting}), $\mu=\xi\left (\gener\right )$ lifts to the element 
  $\sigma =\xi (\beta _{\mu})$ of $\chi (M)$, where
  \begin{displaymath}
  \beta _{\mu}\left (\sum_{k=0}^{2}{a_{k} v^{k}}\right )=\sum_{k=0}^{2}{(\gener)(a_{k})
  v^{k}}.
  \end{displaymath}
  Since the only possibilities for $\chi (M)$ are $\luna$ and $\gruppo\oplus\gruppo$, it is enough to show that 
  $\sigma ^{3}\neq 1$, i.e. $\beta _{\mu} ^{3} \not\in\innerM$.
  
  From the relations
  $$
  (\alpha\otimes\beta )(1\otimes (U_{0}^{*})^{3}g)=\obstruction\, (1
  \otimes (U_{0}^{*})^{3}g)
  $$
  and
  $$
  \Ad W (1\otimes (U_{0}^{*})^{3} g)=1\otimes (U_{0}^{*})^{3}g
  $$ 
  we obtain that 
  \begin{equation*}
  \gamma (1\otimes (U_{0}^{*})^{3}g)= \obstruction (1\otimes (U_{0}^{*})^{3}g). 
  \end{equation*}
  Since $\Ad v=\gamma$\, on $\tenpro$, the last equality can be rewritten as  
  \linebreak $\Ad v(1\otimes(U_{0}^{*})^{3}g)=\obstruction (1\otimes (U_{0}^{*})^{3}\,
  g)$ or
  \begin{equation}
  \label{vact}
  \Ad (1\otimes (U_{0}^{*})^{3}\, g)(v)=\obstruconj v.
  \end{equation}
  
  Thus, using the definition of $\beta _{\mu}$, the relation $\beta ^{3}=\Ad g$,  and (\ref{vact}) we obtain
  \begin{displaymath}
  \begin{split}
  \beta _{\mu} ^{3}\left (\sum_{k=0}^{2}{a_k\, v^{k}}\right )& =\sum_{k=0}^{2}
  {\left (Id\otimes(\Ad (U_{0}^{*})^{3}\, \beta ^{3})\right )(a_{k})v^{k}}=\sum_{k=0}^{2}{\ad (1\otimes (U_{0}^{*})^{3}g)(a_{k}) 
  v^{k}} \\
  &=\Ad (1\otimes (U_{0}^{*})^{3}g)\left (\sum_{k=0}^{2}e^{\frac{2\pi 
  i k}{3}}a_{k}\, v^{k}\right ), 
  \end{split}
  \end{displaymath}
  which implies that up to an inner automorphism $\beta _{\mu}^{3}$ is a dual action, so it is outer and $\chi (\M)\cong\luna$.

\qvd

\section{$\crossprod$ is not anti-isomorphic to itself.}

In this section we are going to show that $\M=\crossprod$ is not 
anti-isomorphic to itself by using the dual action $\dualgr\to\autM$, 
which gives rise to the only subgroup of order $3$ in $\chi (\M)$. This 
argument has been described by Connes in \cite{Connes4} and \cite{Connes6}.

First of all note that the action $\gamma$ can be decomposed as 
$$
\gamma=\ad W\,\gamma_{1}\gamma_{2},
$$ 
where $\gamma_{1}\in\Ctpten$ and
$\gamma_{2}\in\chiusopten$ and $W$ is a unitary in $\tenpro$. 

In fact, $\gamma=\Ad W(\alpha\otimes Id)(Id\otimes\beta)$ and 
$\gamma _{1}=\alpha\otimes Id$ is centrally trivial, since any central 
sequence in $\tenpro$ has the form $(1\otimes x_{n})$, for a central sequence 
$(x_n)$ in $R_{0}$. Furthermore in the proof of Lemma \ref{eleventh} we showed 
that $\beta=\displaystyle\lim_{n\rightarrow\infty}{\Ad x_{n}}$ so that 
$\gamma _{2}=Id\otimes\beta=\lim_{n\rightarrow\infty}{\Ad (1\otimes x_{n})}$ 
belongs to $\chiusopten$. 

Note also that this decomposition of $\gamma$ into an approximately inner 
automorphism and a centrally trivial automorphism is unique up to inner 
automorphisms, since $\chi\pten=1$. 

Let $M$ be an arbitrary von Neumann algebra. Define the conjugate $M^{c}$ 
of $M$ as the algebra whose underlying vector space is the conjugate of $M$
(i.e. for $\lambda\in\mathbb{C}$ and  $x\in M$ the product of $\lambda$ by $x$ in $M^{c}$
is equal to $\bar{\lambda}x$) and whose ring structure is the same as in $M$.
The opposite $M^{o}$ of $M$ is by definition the algebra whose underlying 
vector space is the same as for $M$ while the product of $x$ by $y$ is equal 
to $yx$ instead of $xy$. $M^{c}$ and $M^{o}$ are clearly isomorphic through the
map $x\rightarrow x^{*}$. For $\psi\in\autM$ we denote by $\psi ^{c}$ the 
automorphism of $M^{c}$ induced by $\psi$. 

For the convenience of the reader we detail here Connes's argument to show that the factor $\crossprod$ is not antiisomorphic to itself (see \cite{Connes4} and 
\cite{Connes6} for Connes's  argument).

\begin{theorem}
   $\M=\crossprod$ is not anti--isomorphic to itself.
\end{theorem}
{\bf Proof} 
    We proved in the previous section (Theorem \ref{twelve}) that 
    $\chi (\M)\cong\luna$ and that the dual action $\widehat{\gamma}:\dualgr\to\autM$
    produces the only subgroup of order $3$ in $\chi (\M)$, namely $\langle\sigma 
    ^{3}\rangle =\langle\xi (\beta _{\mu} ^{3})\rangle$.

    Since $\chi (\M)$ is an invariant of the von Neumann algebra $\M$, it follows 
    that $\dualgr$ is an invariant of $\M$. This implies that
    $\dualcr$ is an invariant of $\M$, as it is the dual action $\widetilde{\gamma}:\gruppo\longrightarrow\mbox{Aut}(\dualcr )$ of $\widehat{\gamma}$. 

    Since, by  Takesaki's duality theory [Theorem $4.5$, \cite{Take}] 
    \[ \dualcr \cong\pten\otimes B(\ell ^{2}(\gruppo )), \]
    and in this identification the dual action $\widetilde{\gamma}$ of $\widehat{\gamma}$ 
    corresponds to the action \linebreak $\gamma\otimes\Ad (\lambda (1)^{*})$, where 
    $\lambda$ is the left regular representation of $\gruppo$ on 
    $\ell ^{2}(\gruppo )$, we conclude that $\gamma$ is an invariant of $M$. 
    Consequently, the centrally trivial automorphism $\gamma _{1}$ and the approximately inner automorphism $\gamma_{2}$, which 
    appear in the decomposition of $\gamma$, are also invariants of $\M$. 

    Note that 
    \[ \gamma _{1}^{3}=\Ad (u\otimes 1)\quad\text{and}\quad\gamma_{1}(u\otimes 
    1)=\obstruction\, (u\otimes 1), \]
    so that $\gamma_{1}$ is a $\gruppo$--kernel of $\tenpro$ with obstruction
    $\obstruction$ to lifting.

    Observe also that in the above argument we have only used the abstract
    group $\dualgr$ and the dual action defined on $\dualcr$. Hence we have
    found a canonical way to associate to the von Neumann algebra $\M$ a 
    scalar, equal to $\obstruction$ in our case, which is invariant under 
    isomorphisms.
    
    Now if $\M=\crossprod$ was anti--isomorphic to itself, then $\M$ and 
    $\M ^{c}$ would be isomorphic. But the obstruction 
    associated to $\gamma _{1}^{c}$, and therefore to $\M ^{c}$ is equal to 
    $\obstruconj$.
   
\qvd

\section*{Acknowledgments}

We thank our advisor, Professor F. R\u{a}dulescu for his help and support,
Professor D. Bisch and Professor D. Shlyakhtenko for useful discussions,
and Dr. M. M\"{u}ger for his suggestions regarding this manuscript.

\end{document}